\documentclass[a4paper]{cas-dc}
\usepackage[numbers]{natbib}
\usepackage{graphicx}%
\usepackage{multirow}%
\usepackage{amsmath,amssymb,amsfonts}%
\usepackage{amsthm}%
\usepackage{mathrsfs}%
\usepackage[title]{appendix}%
\usepackage{xcolor}%
\usepackage{textcomp}%
\usepackage{manyfoot}%
\usepackage{booktabs}%
\usepackage{algorithm}%
\usepackage{algorithmicx}%
\usepackage{algpseudocode}%
\usepackage{listings}%
\usepackage{tabularray}
\usepackage{float}
\usepackage{subfigure}
\usepackage[pagewise]{lineno}
\usepackage{hyperref} 
\newtheorem{lemma}{Lemma}
\newtheorem{theorem}{Theorem}
\bibpunct[, ]{[}{]}{,}{n}{,}{,}
\renewcommand\bibfont{\fontsize{10}{12}\selectfont}

\begin{document}

\title{Efficient Computation of the Non-convex Quasi-norm Ball Projection with Iterative Reweighted Approach}
\shorttitle{Efficient Projection onto Non-Convex Quasi-Norm Balls}
\author[1,2]{Qi An}

\author[1]{Jiao Wang}

\author[1]{Zequn Niu}

\author[3]{Nana Zhang}
\cormark[1]

\affiliation[1]{organization={Department of Computer Science, The Open University of China}, addressline={75 Fuxing Road}, city={Beijing}, postcode={100039}, country={China}}

\affiliation[2]{organization={Engineering Research Center of Integration and Application of Digital Learning Technology, Ministry of Education}, addressline={2 Weigongcun Road}, city={Beijing}, postcode={100081}, country={China}}

\affiliation[3]{organization={College of Economics \& Management, Zhejiang University of Water Resources and Electric Power}, addressline={538 Xuelin Road}, city={Hangzhou}, postcode={310018}, country={China}}

\cortext[cor1]{Corresponding author}


\begin{abstract}
In this study, we focus on computing the projection onto the $\ell_p$ quasi-norm ball, which is challenging due to the non-convex and non-Lipschitz nature inherent in the $\ell_p$ quasi-norm with $0<p<1$. We propose a novel localized approximation method that yields a Lipschitz continuous concave surrogate function for the $\ell_p$ quasi-norm with improved approximation quality. Building on this approximation, we enhance the  state-of-the-art iterative reweighted algorithm proposed by Yang et al. (J Mach Learn Res 23:1–31, 2022) by constructing tighter subproblems. This improved algorithm solves the $\ell_p$ quasi-norm ball projection problem through a series of tractable projections onto the weighted $\ell_1$ norm balls. Convergence analyses and numerical studies demonstrate the global convergence and superior computational efficiency of the proposed method. 
\end{abstract}

\maketitle

\section{Introduction}
Sparse learning problems are often cast as non-convex $\ell_p$ optimization problems, leveraging the $\ell_p$ quasi-norm to promote both low-rankness and model sparsity \cite{Chartrand2007, lu2017, 6751136}. Although non-convex optimization problems are generally challenging, fields like machine learning and signal processing often feature problem structures that allow effective solutions to be found. Projected gradient descent (PGD), a cornerstone for solving such structured non-convex constrained optimization, whose practicability and efficiency rely heavily on the projection step. Consequently, recent algorithmic advancements \cite{Drusvyatskiy2015, Hager2014ProjectionAF} have focused on computing Euclidean projections of signals onto non-convex sets. In this paper, we specifically focus on the efficient projections of a signal (vector) in $\mathbb{R}^n$ onto the $\ell_p$ quasi-norm ball with $0<p<1$. This problem can be formulated as follows.

\begin{equation}
    \min\limits_{\boldsymbol{x}\in\mathbb{R}^n}~\frac{1}{2}\Vert \boldsymbol{x}-\boldsymbol{y} \Vert^2_2 \qquad
    \text{s.t.}~\Vert \boldsymbol{x} \Vert^p_p\leq\gamma, \label{eq1}
\end{equation}
where $\boldsymbol{y}\in\mathbb{R}^n$ is the signal to be projected and $\gamma>0$ is a scalar specifying the radius of the $\ell_p$ quasi-norm ball. Choosing a value of $p$ where $0<p<1$ makes $\Vert \boldsymbol{x} \Vert_p=(\sum_{i=1}^n{|x_i|}^p)^{1/p}$ behave as a quasi-norm rather than a proper norm of vector $\boldsymbol{x}$ \cite{chen2011}. To avoid trivial cases, it is assumed that $\Vert \boldsymbol{y} \Vert^p_p>\gamma$. Various efficient methods have recently been proposed to address problem (\ref{eq1}).

%applicability
Computing the projection onto a $\ell_p$ quasi-norm ball is crucial for addressing various inverse problems in fields such as computer vision, signal processing, and machine learning \cite{balda2020adversarial, nikolova2005analysis, thom2015efficient}. These problems often involve minimizing a loss function to evaluate model performance while imposing a complexity constraint on the model vector. Recently, there has been a surge in interest in sparse learning models due to their desirable properties, such as reduced system complexity and improved generalization performance. Typically, these models incorporate a sparsity-inducing constraint or a regularizer, often using the $\ell_1$ norm \cite{Tibshirani2005SparsityAS}. 
However, $\ell_p$ quasi-norm ($0<p<1$), known for its non-convex and non-Lipschitz nature, has gained traction due to its ability to induce superior model sparsity across a wide range of applications, including sparse recovery, matrix completion, and compressed sensing \citep{wang2021nonconvex}. The sparsity-inducing approach is typically formulated as an   $\ell_p$-constrained minimization problem, 
\begin{equation}
 \min\limits_{\boldsymbol{x}\in\mathbb{R}^n}~\psi(\boldsymbol{x}; \boldsymbol{y}) \qquad
    \text{s.t.}~\Vert \boldsymbol{x}\Vert^p_p\leq\gamma, \label{eq2} 
\end{equation}
where $\psi(\boldsymbol{x}; \boldsymbol{y})$ evaluates how well the model vector $\boldsymbol{x}\in\mathbb{R}^n$ fits the observed signal $\boldsymbol{y}\in\mathbb{R}^m$. Recovering inherently sparse vectors is fundamental to many machine learning applications such as signal denoising and data compression \cite{4011956}. Imposing an $\ell_p$ quasi-norm ball constraint is an effective approach
to promote sparsity in the unknown variables. A notable special case of this problem is the classical sparse coding problem, 
\begin{equation}
    \min\limits_{\boldsymbol{x}\in\mathbb{R}^n}~\Vert A\boldsymbol{x}-\boldsymbol{y} \Vert_2^2 \qquad
    \text{s.t.}~\Vert \boldsymbol{x}\Vert^p_p\leq\gamma, \label{eq22} 
\end{equation}
with $A\in\mathbb{R}^{n \times m}$ being an over-complete dictionary.
Using the $\ell_p$ quasi-norm ($0<p<1$) in the sparsity constraint poses significant computational challenges, and the sparsity-constrained least squares problem (\ref{eq22}) is proved to be NP-hard \cite{Natarajan1995}. As a popular workaround, a relaxation-based technique known as the iterative PGD method was introduced by \citep{BAHMANI2013366} to address (\ref{eq22}). The efficacy of the PGD algorithm relies on the existence of practical heuristics for computing the projection onto an $\ell_p$ quasi-norm ball. Consequently, developing an effective algorithm for solving the projection problem (\ref{eq1}) is essential, as it forms the foundation for designing iterative algorithms for solving optimization problems that use the $\ell_p$ quasi-norm to induce sparsity. Moreover, with quasi-norm spherical projections embedded as subprograms within the algorithmic framework, fast convergence and low complexity in these computations are required.

The non-convex and non-Lipschitz nature of $\Vert \boldsymbol{x}\Vert^p_p$ introduces substantial challenges in deriving verifiable optimality conditions and defining optimal solutions, making it difficult to design efficient solving algorithms. Consequently, research efforts \cite{chen2012a, chenng2012} have shifted towards solving an alternative formulation with an $\ell_p$ quasi-norm regularizer, expressed as
\begin{equation}
    \min\limits_{\boldsymbol{x}\in\mathbb{R}^n}~\frac{1}{2}\Vert \boldsymbol{x}-\boldsymbol{y} \Vert^2_2 +
\theta \Vert \boldsymbol{x} \Vert^p_p.
\label{eq3}
\end{equation}
Here, the regularization parameter $\theta$ functions similarly to $\gamma$ in problem (\ref{eq1}), regulating the influence of the quasi-norm-based prior. Several algorithms have been proposed for solving the $\ell_p$ quasi-norm regularization problems   \citep{Chen2010ConvergenceOR, Lai2011AnU, lu2014iterative}. For $p=1$, there is a direct correspondence between $\gamma$, $\theta$, and the optimal solutions for both (\ref{eq1}) and { (\ref{eq3})}, allowing algorithms designed for solving (\ref{eq3}) to address (\ref{eq1}) via bisection over $\theta$. However, for $0<p<1$, such correspondence no longer exists and the non-convexity of the $\ell_p$ quasi-norm presents a more general challenge to the constrained formulation  \citep{Gupta2013nonconvexPP}, necessitating an algorithm specifically tailored to solve problem (\ref{eq1}). 
%yet there are limited algorithms available for this purpose.

%solution
Despite this need, there are limited algorithms available for solving the $\ell_p$ quasi-norm constrained problem (\ref{eq1}). While Bahmani and Raj discussed theoretical aspects of projected points, { the authors} did not provide a specific algorithm \citep{BAHMANI2013366}. Gupta and Kumar introduced an exhaustive search algorithm employing  branch-and-bound techniques, which can yield a globally optimal solution with high probability for most values of $p$ between 0 and 1  \citep{Gupta2013nonconvexPP}. This algorithm branches out based on the smallest nonzero components of the projection vector, solving a set of univariate nonlinear equations to find the optimal solution within each branch. Building on this, Chen et al. used the Lagrangian method within a three-block alternating direction method of multipliers (ADMM) framework, treating the $\ell_p$ quasi-norm ball projection as a subproblem of a low-rank matrix decomposition problem \citep{chen2021}. { The authors} created a bisection method to find the optimal Lagrange multiplier, developed a unified expression for the optimal solution, and eliminated redundant branches through iterative projection and Lagrange multiplier estimation. Won et al. investigated a univariate proximal map related to the Lagrangian dual of the reformulated primal problem and developed a bisection approach that delivers accurate solutions to (\ref{eq1}) with minimal duality gap \citep{won2022unified}. 

In a recent advancement, Yang et al. introduced the innovative iterative reweighted $\ell_1$ norm ball projection (IRBP) algorithm \citep{yang2022towards}. This iterative approach resolves a sequence of smooth and convex subproblems to obtain an approximate solution for problem (\ref{eq1}). Specifically, at the $k$th iteration, with current iterate $\boldsymbol{x}^k$ and perturbation $\epsilon^k$,  this approach approximates $\Vert \boldsymbol{x} \Vert_p^p$ using the smooth function $\sum_{i=1}^n (x_i+\epsilon^k)^p$, and relaxes the smoothed quasi-norm ball to a weighted $\ell_1$ norm ball with respect to $\boldsymbol{x}^k$. The next iterate $\boldsymbol{x}^{k+1}$ is determined as the projection onto this derived weighted $\ell_1$ norm ball. The sequence of perturbation diminishes over iterations, ensuring that the relaxed balls $\ell_1$ norm can approach the original $\ell_p$ quasi-norm ball. The  IRBP method follows the majorization-minimization scheme \citep{lange2016}, a common algorithmic framework for solving non-convex and non-smooth minimization problems that involve $\ell_p$ quasi-norm-based sparsity-inducing terms \citep{Chartrand2008, Figueiredo2007}. This method guarantees convergence to a feasible first-order stationary point, supported by a dynamic updating strategy for the perturbation parameter. Since the weighted $\ell_1$ ball projection can be efficiently solved and is guaranteed to terminate after a finite number of iterations, the IRBP algorithm exhibits superior computational performance compared to other state-of-the-art methods.

In this paper, we propose an enhancement to the existing IRBP algorithm designed for solving the projection onto the $\ell_p$ quasi-norm ball. The key innovation of our approach involves devising a novel method to provide $\epsilon$-approximation $\Vert x \Vert_p^p$ to handle its non-differentiability, which provides a more accurate approximation under the same perturbation magnitude.

Leveraging this novel smoothing technique, we adapt the existing iterative algorithm to solve the projection onto the $\ell_p$ quasi-norm ball more efficiently. Specifically, we construct a sequence of weighed $\ell_1$ norm balls that offer improved approximation of the $\ell_p$ quasi-norm ball. We also establish that any accumulation point of the sequence produced by the modified algorithm remains a first-order stationary solution, provided that a proper strategy is employed for updating the perturbation parameter.

Throughout this paper, we adhere to the following notation conventions. Lowercase letters such as $k$, $n$ indicate scalars, while bold lowercase letters such as $\boldsymbol{x}$,  $\boldsymbol{y}$ represent vectors. The notation $x_i$ refers to the $i$th component of vector $\boldsymbol{x}$.  To simplify notation, the set $\{1, 2, \cdots,n\}$ is denoted by $[n]$. We use $\boldsymbol{x}_{\mathcal{I}}$ to denote the sub-vector of $\boldsymbol{x}$ containing entries indexed by $\mathcal{I}\subset[n]$. The set of zero entries and non-zero entries of vector $\boldsymbol{x}$ are represented by $\mathcal{A}(\boldsymbol{x})$ and $\mathcal{I}(\boldsymbol{x})$, respectively. To generalize, we further define $\mathcal{A}_\epsilon(\boldsymbol{x})=\{i\in[n] \mid x_i\leq\epsilon\}$ and $\mathcal{I}_\epsilon(\boldsymbol{x})=\{i\in[n] \mid x_i>\epsilon\}$. The vector $|\boldsymbol{x}|$ has the same dimension as $\boldsymbol{x}$, with the $i$th component being $|x_i|$. Given a scalar $\tau$, $[\boldsymbol{x}]^\tau$ denotes a vector with the same dimension as $\boldsymbol{x}$, where the $i$th component is $x_i^\tau$. Let $\mathbb{R}^n_+$ be the nonnegative orthant of $\mathbb{R}^n$, and  $\mathbb{R}^n_{++}$ denotes the interior of $\mathbb{R}^n_+$.

The organization of this paper is as follows. Section 2 introduces an alternative  iteratively reweighted $\ell_1$ norm algorithm that achieves a more accurate approximation for solving the projection onto the $\ell_p$ quasi-norm ball. Section 3 analyzes the global convergence properties of the proposed algorithm. Section 4 presents the results of numerical experiment. Section 5 provides concluding remarks.

\section{Projection onto the $\ell_p$ quasi-norm ball}
We first review the properties of the projection of a signal onto an $\ell_p$ ($0<p<1$) quasi-norm ball and the optimality conditions of the projection problem. Regarding problem (\ref{eq1}), Bahmani and Raj identified several properties that the optimal solution must meet, as summarized in the following lemma.\\

\begin{lemma} %lemma1
\cite{BAHMANI2013366} Given$\Vert \boldsymbol{y} \Vert_p^p>\gamma$, let $\boldsymbol{x}^*$ be an optimal solution of problem (\ref{eq1}). The following statements hold:
\renewcommand{\labelenumi}{(\roman{enumi})}
\begin{enumerate}
\item $\Vert \boldsymbol{x}^* \Vert_p^p=\gamma$.
\item $x_i^*y_i\geq0$ and $|x_i^*|\leq |y_i|$ for all $i\in[n]$.
\item If $|y_i|>|y_j|$ for some $i,~j\in[n]$, then $|x_i^*|>|x_j^*|$.
\end{enumerate}
\label{lemma1}
\end{lemma}

It follows from Lemma \ref{lemma1} that the projection $\boldsymbol{x}^*$ shares the same support as the signal $\boldsymbol{y}$, with the signs of individual components matching those of the corresponding components of $\boldsymbol{y}$. Thus, when handling a point $\boldsymbol{y}$ containing positive, zero, and negative components, we can first compute the projection of $|\boldsymbol{y}|_{\mathcal{I}(\boldsymbol{y})}$ and later incorporate the zero components and signs. Consequently, without loss of generality, we assume that $\boldsymbol{y}\in\mathbb{R}^n_{++}$. This assumption allows us to simplify the optimization problem (\ref{eq1}) as

\begin{equation}
    \min\limits_{\boldsymbol{x}\in\mathbb{R}^n}~\frac{1}{2}\Vert \boldsymbol{x}-\boldsymbol{y} \Vert^2_2 \qquad
\text{s.t.}~\Vert \boldsymbol{x} \Vert^p_p\leq\gamma,~\boldsymbol{x}\geq0. 
\tag{P}
\label{eq4}
\end{equation}

Any point $\boldsymbol{x}$ satisfies the first-order optimality conditions of problem (\ref{eq4}) if and only if 
\begin{align*}
    (y_i-x_i)x_i&=\lambda px_i^p,\quad i\in[n],\\
    \sum_{i=1}^n x_i^p&=\gamma,\\
    \boldsymbol{x}\geq\boldsymbol{0},\quad &\lambda \geq 0,
\end{align*}
where $\lambda$ represents the dual variable.

{Our investigation begins by introducing a more sophisticated smoothing technique to handle $\Vert\boldsymbol{x}\Vert_p^p$. It is well-established that for $0<p<1$, $|t|^p$ is not locally Lipschitz continuous at 0, although it is locally Lipschitz continuous at all other points. This poses a significant challenge in designing algorithms to solve problem (\ref{eq4}), as Lipschitz continuity is a necessary condition for obtaining the Clarke subgradient. To address this issue, several approximation methods have been explored. For instance, existing IRBP methods \cite{yang2022towards, An09122024} replaced $\Vert\boldsymbol{x}\Vert_p^p$ with smoothing functions such as $\sum_{i=1}^n (|x_i|+\epsilon)^{p}$ or $\sum_{i=1}^n (|x_i|^2+\epsilon^2)^{p/2}$, where $\epsilon$ is a small perturbation parameter. These approximations ensure Lipschitz continuity by smoothing the function across the entire domain. However, since the non-Lipschitz behavior occurs only at 0, we propose, in this work, a localized approximation approach that targets only the neighborhood around 0, while ensuring Lipschitz continuity is preserved throughout the entire domain.}

With this motivation, we turn our attention to the following approximation to $|t|^p$ proposed by \cite{lu2014iterative}:

\begin{equation}
    h_u(t):=\min_{0\leq s\leq u} p(|t|s-\frac{s^q}{q}), \quad\forall t \in \mathbb{R}.
    \label{luequation}
\end{equation}
where $u>0$ is arbitrarily chosen and $q$ satisfies $\frac{1}{p}+\frac{1}{q}=1$. Let $g_t(s):= p(|t|s-\frac{s^q}{q})$ for $s>0$. Since $g_t$ is continuous and $g_t(s)\rightarrow\infty$ as $s\downarrow0$, the function $h_u(t)$ is well-defined for all $t\in\mathbb{R}$. Lu introduced this $\epsilon$-approximation in the development of an iterative reweighted minimization method for solving $\ell_p$ quasi-norm regularized unconstrained problems \citep{lu2014iterative}. Xiu et al. further generalized this iterative reweighted method to solve the $\ell_1-\ell_p$ minimization problems \cite{xiu2018}.
The Lipschitz continuity of $h_u$ and the approximate error between $|t|^p$ and $h_u(t)$ were also provided, as demonstrated in the following lemma.\\

\begin{lemma}\cite{lu2014iterative}
Let $u>0$ be arbitrarily given, the following statements hold:
\renewcommand{\labelenumi}{(\roman{enumi})}
\begin{enumerate}
\item $h_u$ is $pu$-Lipschitz continuous on $(-\infty,\infty)$.
\item $0\leq h_u(t)-|t|^p \leq { u^q}$ for every $t\in(-\infty,\infty)$. 
\end{enumerate}
\label{lemma2}
\end{lemma}

Identifying the optimal value of (\ref{luequation}) reveals that $h_u(t)=\phi(|t|;u)$ for all $t$, with $\phi:[0,\infty)\rightarrow(-\infty,\infty)$ defined as
\begin{equation}
    \phi(t;u)=
    \begin{cases}
        t^p & \text{if $t>u^{q-1}$,} \\
        p(tu-u^q/q) & \text{if $0\leq t\leq u^{q-1}$}. \\
    \end{cases}
    \label{eq7}
\end{equation}

The following lemma establishes several proprieties of the differentiable function $\phi$.\\

\begin{lemma}
The following statements hold:
\renewcommand{\labelenumi}{(\roman{enumi})}
\begin{enumerate}
\item Let $u>0$ be arbitrarily given, $\phi(t;u)$ is concave on $[0,\infty)$. 
\item Let $u_1>u_2>0$ be arbitrarily given, $\phi(t;u_1)\leq\phi(t;u_2)$ for all $t\in[0,\infty)$.
\end{enumerate}
\label{lemma3}
\end{lemma}
\begin{proof}
(i) We establish the concavity of $\phi$  by confirming that its first derivative is non-increasing. This is evident upon characterizing the first derivative of $\phi$ as 
\begin{equation}
    \phi'(t;u)=
    \begin{cases}
        pt^{p-1} & \text{if $t>u^{q-1}$} \\
        pu & \text{if $0\leq t\leq u^{q-1}$}\\
    \end{cases},
    \label{eq8}
\end{equation}\\
and recognizing the relationship $(q-1)(p-1)=1$.

(ii) When $u_1>u_2$, it follows that $u_1^{q-1}<u_2^{q-1}$. By definition,  $\phi(t;u_1)=\phi(t;u_2)$ holds for all $t\in[u_2^{q-1},+\infty)$. For all $t\in[0,u_2^{q-1})$, we have $\phi'(t;u_1)> \phi'(u_2^{q-1};u_1) = pu_2 = \phi'(t;u_2)$. This, combined with $\phi(u_2^{q-1};u_1)=\phi(u_2^{q-1};u_2)$, implies that $\phi(t;u_1)<\phi(t;u_2)$ for all $t\in[0,u_2^{q-1})$.
\end{proof}

By the Lipschitz continuity of $h_u$ over $(-\infty, \infty)$, substituting $|x_i|^p$ with $h_u(x_i)$ for all $i\in[n]$ yields a Lipschitz continuous $\epsilon$-approximation to the $\ell_p$ quasi-norm of any $\boldsymbol{x}\in\mathbb{R}^n$:
\begin{equation}
\Vert \boldsymbol{x} \Vert^p_p \approx \sum_{i=1}^n h_u(x_i).
\end{equation}

{For all values of $u$, $h_u(x_i)$ provides a local approximation to $|x_i|^p$ within a neighborhood around 0. Rather than approximating the function globally, it approximates $|x_i|^p$ only when $|x_i|$ is below a threshold value $u^{q-1}$, while Lipschitz continuity is preserved in the entire domain. Furthermore, we observe that by selecting an appropriate value of $u$, we can leverage $h_u$ to achieve a more refined approximation of the $\ell_p$ quasi-norm.} Specifically, for any perturbation parameter $\epsilon>0$, we define $u_\epsilon:=\epsilon^{p-1}$, and the smoothing function $h_{u_\epsilon}(x_i)$ offers an enhanced local approximation to $|x_i|^p$ compared to the commonly used approximation $(x_i+\epsilon)^p$, as illustrated in the next lemma. \\

\begin{lemma}%%lemma4
Let $\epsilon>0$ be arbitrarily given and $u_\epsilon=\epsilon^{p-1}$. For all $t \in [0,\infty)$, it holds that $\phi(t;u_\epsilon) \leq (t+\epsilon)^{p}$. 
%In particular, for all $t \in [0, u_\epsilon^{q-1})$, $\phi(t;u_\epsilon) \leq (1-p)(t+\epsilon)^{p}$ holds.
\label{lemma4}
\end{lemma}

\begin{proof}
For any $u>0$, since $\phi'(t;u)-p(t+\epsilon)^{p-1}$ is increasing over the interval $[0,u^{q-1})$, $\phi(t;u)-(t+\epsilon)^{p}$ is convex on $[0,u^{q-1}]$ due to the increasing slope. With respect to $u_\epsilon=\epsilon^{p-1}$, we have $\phi(0;u_\epsilon)=(1-p)\epsilon^{p}<\epsilon^{p}$ and $\phi(u_\epsilon^{q-1};u_\epsilon)=\epsilon^p< (u_\epsilon^{q-1}+\epsilon)^{p}$. By convexity, it can be inferred that $\phi(t;u_\epsilon)<(t+\epsilon)^{p}$ holds over the interval $[0,u_\epsilon^{q-1})$. Moreover, for all $t\in[u_\epsilon^{q-1}, +\infty]$, $\phi(t;u_\epsilon) \leq (t+\epsilon)^{p}$ directly follows from the definition. 
%Given $\phi(0;u_\epsilon)=(1-p)\epsilon^{p}$, to establish $\phi(t;u_\epsilon) \leq (1-p)(t+\epsilon)^{p}$ for all $t \in [0, u_\epsilon^{q-1})$, it suffices to demonstrate $\phi(t;u_\epsilon) - (1-p)(t+\epsilon)^{p}$ is increasing over  [0, u_\epsilon^{q-1}), or,  $\phi'(t;u_\epsilon)-p(1-p)(t+\epsilon)^{p-1}>0$, for all $t \in [0, u_\epsilon^{q-1})$. This follows from  $pu_\epsilon-p(1-p)(t+\epsilon)^{p-1}>pu_\epsilon-p(1-p)\epsilon^{p-1}>0$. 
\end{proof}

From Lemma \ref{lemma2} and Lemma \ref{lemma4}, it can be inferred that with the same amount of perturbation $\epsilon$, substituting $h_{u_\epsilon}$ for $|\cdot|^p$ yields a more effective approximation for $\Vert \boldsymbol{x} \Vert^p_p$. In other words, we have the following inequality:
\begin{equation}
    \sum_{i=1}^n |x_i|^{p} \leq 
    \sum_{i=1}^n h_{u_\epsilon}(x_i) \leq
    \sum_{i=1}^n (|x_i|+\epsilon)^{p}.
    \label{eqsum}
\end{equation}

Building on the smoothing function $h_{u_\epsilon}$ that provides refined local approximation, we proceed to formulate iterative subproblems as a sequence of weighted $\ell_1$ norm ball projection problems, following the approach of the IRBP algorithm proposed in \cite{yang2022towards}. In these subproblems, $\Vert \boldsymbol{x} \Vert_p^p$ is approximated by $\sum_{i=1}^n h_{u_{\epsilon^k}}(x_i)$, where $\{\epsilon^k\}$ is a sequence of positive values gradually converging to zero as $k \rightarrow \infty$. 

At the $k$th iteration with perturbation $\epsilon^k>0$, when employing the smoothing function $\sum_{i=1}^n  h_{u_{\epsilon^k}}(x_i)$ to approximate $\Vert \boldsymbol{x} \Vert_p^p$ in the original problem (\ref{eq4}), we encounter a projection problem essentially taking the form
\begin{equation}
    \min\limits_{\boldsymbol{x}\in\mathbb{R}^n}~\frac{1}{2}\Vert \boldsymbol{x}-\boldsymbol{y} \Vert^2_2 \qquad
    \text{s.t.}~ \sum_{i=1}^n \phi^k(x_i)\leq \gamma,~\boldsymbol{x}\in\mathbb{R}^n_+, \label{eq9}
\end{equation}
where $\phi^k(\cdot)$ is an abbreviation for $\phi^k(\cdot; u_{\epsilon^k})$.

%The smoothed constraint in problem (\ref{eq9}) can be linearized without violating the feasibility of the original problem (\ref{eq4}). 
Leveraging the concavity of $\phi^k$, we derive an upper bound on $\phi^k(x_i)$ for all $x_i\in[0,+\infty)$ with respect to the current iterate $\boldsymbol{x}^k$, which leads to
\begin{equation}
\phi^k(x_i) \leq \phi^k(x_i^k)+(\phi^{k})'(x_i^k)(x_i-x_i^k).
\label{eq10}
\end{equation}
Substituting this upper bound and $(\phi^{k})'(x_i^k)$, the smoothed constraint in problem (\ref{eq9}) can be relaxed to \begin{equation*}
\sum_{i=1}^n \phi^k(x_i^k)+ \sum_{i\in\mathcal{I}_{\epsilon^k}(\boldsymbol{x}^k)} 
p{(x_i^k)}^{p-1}(x_i-x_i^k) +\sum_{i\in\mathcal{A}_{\epsilon^k}(\boldsymbol{x}^k)} 
p{(\epsilon^k)}^{p-1}(x_i-x_i^k) \leq\gamma,
\end{equation*}
which constitutes a localized linear approximation of the original $\ell_p$ quasi-norm constraint. Here, we emphasize that the concavity of $\phi^k$ allows for the linearization of the smoothed constraint, while the refined approximation provided by $\phi^k$ enables the construction of a tighter linear constraint. To simplify, we will denote $\mathcal{I}_{\epsilon^k}(\boldsymbol{x}^k)$ and $\mathcal{A}_{\epsilon^k}(\boldsymbol{x}^k)$ as $\mathcal{I}_\epsilon^k$ and $\mathcal{A}_{\epsilon}^k$, respectively. Defining 
\begin{equation*}
\gamma^k:=\gamma -\sum_{i=1}^n\phi^k(x_i^k)
+\sum_{i\in\mathcal{I}_\epsilon^k} 
p{(x_i^k)}^{p-1}x_i^k +\sum_{i\in\mathcal{A}_\epsilon^k} 
p{(\epsilon^k)}^{p-1}x_i^k,
\end{equation*}
we derive the $k$th convex subproblem denoted as 
\begin{align}
\begin{split}
\min\limits_{\boldsymbol{x}\in\mathbb{R}^n_+}~&\frac{1}{2}\Vert \boldsymbol{x}-\boldsymbol{y} \Vert^2_2 \qquad\\
\text{s.t.}~ \sum_{i\in\mathcal{I}_\epsilon^k} 
p{(x_i^k)}^{p-1}x_i +&\sum_{i\in\mathcal{A}_\epsilon^k} 
p{(\epsilon^k)}^{p-1}x_i \leq \gamma^k. 
\end{split}
\tag{P$_k$}
\label{Pk}
\end{align}
During the iteration process, the new iterate $\boldsymbol{x}^{k+1}$ is obtained as the optimal solution of (\ref{Pk}). The following lemma outlines the relevant properties of $\boldsymbol{x}^{k+1}$.\\

\begin{lemma}%lemma 5
For any $k\in\mathbb{N}$, at the $k$th iteration with $\boldsymbol{x}^k$ and $\epsilon^k>0$, $\boldsymbol{x}^{k+1}$ is the optimal solution of (\ref{Pk}). Provided that $\boldsymbol{x}^k$ is feasible to subproblem (\ref{Pk}), i.e., $\sum_{i=1}^n\phi^k(x_i^k)\leq\gamma$, the following statements hold:
\renewcommand{\labelenumi}{(\roman{enumi})}
\begin{enumerate}
\item $\mathcal{I}(\boldsymbol{x}^{k+1})\neq\varnothing$.
\item $\sum_{i\in\mathcal{I}_\epsilon^k} 
p{(x_i^k)}^{p-1}x_i^{k+1} +\sum_{i\in\mathcal{A}_\epsilon^k} 
p{(\epsilon^k)}^{p-1}x_i^{k+1}=\gamma^k$.
\item $0\leq x_i^{k+1} \leq y_i$ for all $i\in[n]$.
\end{enumerate}
\label{lemma5}
\end{lemma}

\begin{proof}
For the newly constructed subproblem (\ref{Pk}), we obtain that $\gamma^k>\gamma-\sum_{i=1}^n\phi^k(x_i^k)\geq0$, where the last inequality holds by assumption. The remainder of the proof follows directly from the proof of Lemma 7 in \cite{yang2022towards}.

%Since $\gamma^k>0$, subproblem (\ref{Pk}) computes the projection of $\boldsymbol{y}$ onto a weighted $\ell_1$ ball with a non-empty interior. From (\ref{eq10}) we know that the weighted $\ell_1$ ball in the constraint of subproblem (\ref{Pk}) is contained within the perturbed $\ell_p$ norm ball and, therefore, is also contained within the original $\ell_p$ quasi-norm ball. As a result, $\boldsymbol{0}$ cannot be the projection of $\boldsymbol{y}$ onto the weighted $\ell_1$ ball, implying $\mathcal{I}(\boldsymbol{x}^{k+1})\neq\varnothing$. Since $\boldsymbol{y}$ lies outside of the weighted $\ell_1$ norm ball, we can easily prove  by contradiction that $\boldsymbol{x}^{k+1}$, as the projection of $\boldsymbol{y}$ onto the weighted $\ell_1$ norm ball, satisfies $\sum_{i=1}^n wx_i^{k+1}=\gamma^k$ and $0\leq x_i^{k+1} \leq y_i$ for any $i\in[n]$.
\end{proof}

%perturbation $\epsilon^k>0$ satisfying $\phi^k(0)<\gamma/n$
%When $\boldsymbol{x}^k=\boldsymbol{0}$, $\gamma^k=\gamma-\sum_{i=1}^n\phi^k(x_i^k)+\sum_{i=1}^n wi^kx_i^k=\gamma-\sum_{i=1}^n\phi^k(0)>0$. 

Now, we propose a variation of the iterative reweighted $\ell_1$ ball projection algorithm for computing the Euclidean projection of a given signal onto the $\ell_p$ quasi-norm ball. 

\begin{algorithm}
\caption{  }
\label{lloyd} 
\begin{algorithmic}[1]
\item[\textbf{Input.}] Parameters $p$, $n$, $\gamma$. Signal vector to be projected $\boldsymbol{y}\in\mathbb{R}_{++}^n$ satisfying $\Vert \boldsymbol{y} \Vert_p^p>\gamma$.
\item[\textbf{Initialization.}] Generate an initial perturbation parameter $\epsilon^0>0$ satisfying  $\phi^0(0)<\gamma/n$. Based on $\epsilon^0$, generate an initial point $\boldsymbol{x}^0\in\mathbb{R}_{+}^n$ satisfying $\sum_{i=1}^n \phi^0(x_i^0)<\gamma$. Choose parameters $\delta\in(0, 1)$, $\tau>1$, and $M>0$. Set $k:= 0$.
\item[\textbf{Step 1.}] Compute $\gamma^k$ and solve subproblem (\ref{Pk}) for the optimal solution $\boldsymbol{x}^{k+1}$.
\item[\textbf{Step 2.}] Check the stopping criteria: if met, the algorithm stops. Otherwise, check if the condition 
\begin{equation}
\Vert \boldsymbol{x}^{k+1}-\boldsymbol{x}^k \Vert_2 \Vert p(\epsilon^k)^{p-1}\text{sign}(\boldsymbol{x}^{k+1}-\boldsymbol{x}^k)\Vert_2^{\tau} \leq M
\label{eqtriger}
\end{equation} holds. If true, set the new perturbation parameter $\epsilon^{k+1}\in(0, \delta\epsilon^{k}]$. Otherwise, set $\epsilon^{k+1}=\epsilon^{k}$. 
\item[\textbf{Step 3.}] Set $k:=k+1$, and return to Step 1. 
\end{algorithmic}
\end{algorithm}

Algorithm \ref{lloyd} continues to use the dynamic updating strategy proposed by \cite{yang2022towards} to determine when to decrease the perturbation parameter. { The $k$th subproblem (\ref{Pk}) dictates the computational complexity per iteration and can be solved using the exact projection algorithms \citep{perez2022efficient, Konstantinos2010} with either $O(n\log n)$ or $O(n)$ operations.}

\section{Convergence analysis}

In this section, we establish the global convergence of the proposed algorithm and demonstrate that any cluster point of the iterative sequence produced by Algorithm \ref{lloyd} is first-order stationary for problem (\ref{eq4}). We clarify here that the convergence analysis in this manuscript parallels the foundational results established in the IRBP framework by \cite{yang2022towards}. Given the modifications introduced by the new smoothing function, however, direct application of their proofs was not feasible. As such, we reconstructed these proofs to align with the properties of Algorithm \ref{lloyd}.

First, we establish that the assumption in Lemma \ref{lemma5} ensuring the feasibility of $\boldsymbol{x}^k$ for subproblem (\ref{Pk}) is consistently maintained. { This, in turn, ensures that $\boldsymbol{x}^k$ remains feasible for the original problem (\ref{eq4}), due to the established inequality (\ref{eqsum}).} \\

\begin{theorem}%%theorem 1
Let $\{\boldsymbol{x}^k\}$ be the sequence produced by Algorithm \ref{lloyd}. For any $k\in\mathbb{N}$, $\boldsymbol{x}^k$ is feasible to subproblem (\ref{Pk}), i.e., $\sum_{i=1}^n \phi^k(x_i^k)\leq\gamma$. 
\label{theorem1}
\end{theorem}
\begin{proof}
We apply induction to prove this. Since $\boldsymbol{x}^0$ is generated such that $\sum_{i=1}^n \phi^0(x_i^0)<\gamma$, the statement is true for $k=0$. Assume that  $\sum_{i=1}^n \phi^k(x_i^k)\leq\gamma$, implying the $k$th subproblem (\ref{Pk}) is feasible since $\boldsymbol{x}^k$ itself is a feasible solution of (\ref{Pk}). Consequently, $\boldsymbol{x}^{k+1}$, as the optimal solution to (\ref{Pk}), satisfies $\sum_{i=1}^n \phi^{k+1}(x_i^{k+1}) \leq \sum_{i=1}^n \phi^k(x_i^{k+1}) \leq \sum_{i=1}^n \phi^k(x_i^k)+(\phi^{k})'(x_i^k)(x_i^{k+1}-x_i^k) \leq \gamma$, where the first inequality is derived from Lemma \ref{lemma3} and the relation $u_{\epsilon^{k}}\leq u_{\epsilon^{k+1}}$, and the second inequality is established through the concavity of $\phi^k$. 
\end{proof}

Theorem \ref{theorem1} confirms that all the subproblems (\ref{Pk}) produced by Algorithm \ref{lloyd} remain feasible, thus ensuring Algorithm \ref{lloyd} is well-defined. Therefore, $\boldsymbol{x}^k$ always stays within the $\ell_p$ { quasi-norm} ball and therefore is feasible to the original problem (\ref{eq4}), according to (\ref{eqsum}). 

Then we establish that the objective of problem (\ref{eq4}) can achieve a quadratic improvement after each iteration. Moreover, irrespective of the strategy for updating $\epsilon^k$, the iterative sequence produced by the proposed algorithm always converges towards a vanishing change in the limit.\\

\begin{lemma}%%lemma 6
    Let $\{\boldsymbol{x}^k\}$ be the sequence produced by Algorithm \ref{lloyd}. It holds that 
    \renewcommand{\labelenumi}{(\roman{enumi})}
    \begin{enumerate}
    \item The sequence $\{\Vert\boldsymbol{x}^k-\boldsymbol{y}\Vert_2^2\}$ monotonically decreases. Indeed, for all $k\in\mathbb{N}$,
    \begin{equation}
        \Vert\boldsymbol{x}^k-\boldsymbol{y}\Vert{}_2^2 - 
        \Vert\boldsymbol{x}^{k+1}-\boldsymbol{y}\Vert_2^2\geq
        \Vert \boldsymbol{x}^k-\boldsymbol{x}^{k+1}\Vert_2^2.
    \end{equation}
    \item The sequence $\{\boldsymbol{x}^k\}$ satisfies $\sum_{k=0}^{+\infty}\Vert \boldsymbol{x}^k-\boldsymbol{x}^{k+1}\Vert_2^2<+\infty$, and consequently, $\lim\limits_{k\rightarrow +\infty} \Vert \boldsymbol{x}^k-\boldsymbol{x}^{k+1}\Vert_2^2=0$.
    \end{enumerate}
    \label{lemma6}
\end{lemma}

\begin{proof}
Since (\ref{Pk}) essentially involves solving the projection onto a closed convex set, the proof follows the proof of Lemma 6 in \cite{yang2022towards}. 
\end{proof} 

To ensure the iterates converge to the first-order stationary points of problem (\ref{eq4}), $\epsilon^k$ must be systematically reduced down to zero throughout the algorithm's execution. This allows $\{\boldsymbol{x}^k\}$ to approach the boundary of the original $\ell_p$ quasi-norm ball in the limit. Therefore, we must establish that the conditon (\ref{eqtriger}) is activated infinitely many times. We denote the set of iterations in which $\epsilon^k$ is decreased as
\begin{equation*}
\mathcal{U}:=\{k\in\mathbb{N}\mid \Vert \boldsymbol{x}^{k+1}-\boldsymbol{x}^k \Vert_2 \Vert p(\epsilon^k)^{p-1}\text{sign}(\boldsymbol{x}^{k+1}-\boldsymbol{x}^k)\Vert_2^{\tau}\leq M\}.
\end{equation*}
Then we have the following result.\\

\begin{lemma} %%lemma 7
It holds that $|\mathcal{U}|= +\infty$, implying $\lim\limits_{k\rightarrow +\infty} \epsilon^k=0$.
    \label{lemma7}
\end{lemma}

\begin{proof}
Assume, for contradiction, that there exists a specific value $\bar{k}\in\mathbb{N}$ such that for all $k\geq \bar{k}$,
\begin{equation}
\Vert \boldsymbol{x}^{k+1}-\boldsymbol{x}^k \Vert_2 \Vert p(\epsilon^k)^{p-1}\text{sign}(\boldsymbol{x}^{k+1}-\boldsymbol{x}^k)\Vert_2^{\tau}>M,
\end{equation}
which implies $\boldsymbol{x}^{k}\neq\boldsymbol{x}^{k+1}$ and $\epsilon^k$ is never reduced for all $k\geq \bar{k}$. Therefore, 
$\Vert \boldsymbol{x}^{k+1}-\boldsymbol{x}^k \Vert_2>M\Vert p(\epsilon^k)^{p-1}\text{sign}(\boldsymbol{x}^{k+1}-\boldsymbol{x}^k)\Vert_2^{-\tau}\geq M(\sum_{x_i^{k+1}\neq x_i^k}p^2{(\epsilon^{\bar{k}})}{}^{2(p-1)}){}^{-\tau/2}>0$ for all $k\geq \bar{k}$, which contradicts Lemma \ref{lemma6}. 
\end{proof}

To simplify, we will denote $\mathcal{I}(\boldsymbol{x}^k)$ and $\mathcal{A}(\boldsymbol{x}^k)$ as $\mathcal{I}^k$ and $\mathcal{A}^k$, respectively. The following lemma provides a preliminary result necessary to establish the proof of global convergence.\\

\begin{lemma} %%lemma 8
Assume 0 is not the cluster point of $\{\lambda^k\}_\mathcal{U}$. Let $\{\boldsymbol{x}^k\}$ be the sequence produced by Algorithm \ref{lloyd}. There exists $\sigma_i>0$ such that $x_i^k>\sigma_i$ for all sufficiently large $k\in\mathcal{U}$ and $i\in\mathcal{I}^k$.
\label{lemma8}
\end{lemma}
\begin{proof}
Since the cluster point of $\{\lambda^k\}_\mathcal{U}$ cannot be 0, there exists $\hat{\lambda}$ such that $\lambda>\hat{\lambda}$ for all $k\in\mathcal{U}$. We prove this by contradiction and suppose this is not true. By the assumption, there exists $\hat{k}$ such that $x_i^k<\frac{1}{2}({\frac{1}{\hat{\lambda} p}\Vert\boldsymbol{y}\Vert_{\infty}})^{1-p}$ and $\epsilon^k<\frac{1}{2}({\frac{1}{\hat{\lambda} p}\Vert\boldsymbol{y}\Vert_{\infty}})^{1-p}$ for all $k>\hat{k}$ and $i\in\mathcal{I}^k$. By Lemma \ref{lemma6} and Lemma \ref{lemma7}, there exists $\bar{k}$ such that $|x_i^k-x_i^{k-1}|<\frac{1}{2}({\frac{1}{\hat{\lambda} p}\Vert\boldsymbol{y}\Vert_{\infty}})^{1-p}$ and $|\epsilon^k-\epsilon^{k-1}|<\frac{1}{2}({\frac{1}{\hat{\lambda} p}\Vert\boldsymbol{y}\Vert_{\infty}})^{1-p}$ for all $k>\bar{k}$ and $i\in\mathcal{I}^k$. Hence for all $k>\max(\hat{k},\bar{k})$ and $i\in\mathcal{I}_{\epsilon^{k-1}}(\boldsymbol{x}^{k-1})$, we have $y_i-x_i^k=\lambda^kp{(x_i^{k-1})}{}^{p-1}\geq \hat{\lambda}p{(x_i^{k-1})}{}^{p-1}\geq \hat{\lambda}p{(x_i^{k}+|x_i^{k-1}-x_i^{k}|)}{}^{p-1}>\Vert\boldsymbol{y}\Vert_{\infty}$; otherwise, we have $y_i-x_i^k=\lambda^kp{(\epsilon^{k-1})}{}^{p-1}\geq \hat{\lambda}p{(\epsilon^{k-1})}{}^{p-1}\geq \hat{\lambda}p{(\epsilon^{k}+|\epsilon^{k-1}-\epsilon^{k}|)}{}^{p-1}>\Vert\boldsymbol{y}\Vert_{\infty}$, which contradicts Lemma \ref{lemma1}.
\end{proof} 

For a primal-dual solution $(\boldsymbol{x}, \lambda)\in\mathbb{R}_+^n\times\mathbb{R}$, the following
metrics are defined to evaluate the optimality residuals at $(\boldsymbol{x}, \lambda)$ for the original problem (\ref{eq4}):
\begin{equation}
    \alpha(\boldsymbol{x}, \lambda):=\sum_{i=1}^n|(y_i-x_i)x_i-\lambda px_i^p|,\quad
    \beta(\boldsymbol{x}):=|\sum_{i=1}^n x_i^p-\gamma|.
\end{equation}
Note that $(\boldsymbol{x}, \lambda)$ is a first-order stationary solution of problem (\ref{eq4}) if and only if $\alpha(\boldsymbol{x}, \lambda)=0$ and $\beta(\boldsymbol{x})=0$. The following lemma establishes the upper bounds for the optimality residuals $\alpha(\boldsymbol{x}^k, \lambda^k)$ and  $\beta(\boldsymbol{x}^k)$.\\

\begin{lemma}%%lemma9
Let $\epsilon^k$ be arbitrarily given, and let $\{\boldsymbol{x}^k\}$ be the sequence produced by Algorithm \ref{lloyd}. Define
    \begin{gather*}
    \mathcal{I}1^k := \{i\in\mathcal{I}^k\mid {x_i^{k-1}\leq \epsilon^{k-1}}\},\\
    \mathcal{I}2^k := \{i\in\mathcal{I}^k\mid \epsilon^{k-1}<x_i^{k-1}\},\\
    \mathcal{I}3^k := \{i\in\mathcal{I}^k\mid {\phi^{k-1}(x_i^{k-1})}>{(x_i^k)}^{p}\},\\
    \mathcal{I}4^k := \{i\in\mathcal{I}^k\mid{\phi^{k-1}(x_i^{k-1})} \leq {(x_i^k)}^{p}\}.
    \end{gather*}
    \renewcommand{\labelenumi}{(\roman{enumi})}
    \begin{enumerate}
    \item There exists $\tilde{\boldsymbol{x}}$ with $\tilde{x}_i^k$ between $x_i^k$ and $\epsilon^{k-1}$ for all $i\in\mathcal{I}1^k$ and $\tilde{x}_i^k$ between $x_i^k$ and $x_i^{k-1}$ for all $i\in\mathcal{I}2^k$ such that 
    $\alpha(\boldsymbol{x}^k, \lambda^k)$ is bounded above by $\lambda^kp(1-p)\epsilon^{k-1}\Vert{[\tilde{\boldsymbol{x}}^k_{\mathcal{I}1^k}]}^{p-2}\Vert_1\Vert \boldsymbol{y}\Vert_\infty+\lambda^kp(1-p)\Vert{[\tilde{\boldsymbol{x}}^k_{\mathcal{I}^k}]}^{p-2}\Vert_2\Vert{\boldsymbol{x}}^{k-1}_{\mathcal{I}^k}-{\boldsymbol{x}}^k_{\mathcal{I}^k}\Vert_2\Vert \boldsymbol{y}\Vert_\infty$.
    \item There exists $\tilde{\boldsymbol{x}}$ with $\tilde{x}_i^k$ between $x_i^k$ and $x_i^{k-1}+\epsilon^{k-1}$ for all $i\in\mathcal{I}3^k$ and $\tilde{x}_i^k$ between $x_i^k$ and $x_i^{k-1}$ for all $i\in\mathcal{I}4^k$ such that 
    $\beta(\boldsymbol{x}^k)$ is bounded above by $p({\epsilon^{k-1}}){}^{p-1}\Vert \boldsymbol{x}^k- {\boldsymbol{x}}^{k-1}\Vert_1 + p\epsilon^{k-1}\Vert{[\tilde{\boldsymbol{x}}_{\mathcal{I}3^k}^k]}^{p-1}\Vert_1 + p\Vert {[\tilde{\boldsymbol{x}}_{\mathcal{I}^k}^k]}^{p-1}\Vert_2\Vert{\boldsymbol{x}}_{\mathcal{I}^k}^{k-1}-\boldsymbol{x}_{\mathcal{I}^k}^k\Vert_2
    +({\epsilon^{k-1}}){}^{p-1} \Vert {{\boldsymbol{x}}^{k-1}_{\mathcal{A}^k}}\Vert_1 + ({\epsilon^{k-1}})^{p}$.
    \end{enumerate}
    \label{lemma9}
    \end{lemma}
    
\begin{proof}
    Since $\boldsymbol{x}^k$ is the optimal solution of (P$_{k-1}$), by the KKT conditons $\boldsymbol{x}^k$ satisfies $y_i-x_i^k-\lambda^k p {(\epsilon^{k-1})}^{p-1} =0$ for all $i\in\mathcal{I}1^k$ and $y_i-x_i^k-\lambda^k p {(x_i^{k-1})}^{p-1} =0$ for all $i\in\mathcal{I}2^k$.     
    We have 
    \begin{align*}
    \alpha(\boldsymbol{x}^k, \lambda^k) 
     =&\sum_{i=1}^n|(y_i-x_i^k)x_i^k-\lambda^k p{(x_i^k)}^{p}|\\
     =&\lambda^kp \sum_{i\in\mathcal{I}1^k} |{(\epsilon^{k-1})}^{p-1}-{(x_i^k)}^{p-1}||x_i^k|\\
      &+\lambda^kp\sum_{i\in\mathcal{I}2^k} |{(x_i^{k-1})}^{p-1}-{(x_i^k)}^{p-1}||x_i^k|.
    \end{align*}
    For any $i\in\mathcal{I}1^k$, by the Lagrange's mean value theorem, there exists $\tilde{x}_i^k$ between $x_i^k$ and $\epsilon^{k-1}$ such that
    \begin{align*}
        |{(\epsilon^{k-1})}^{p-1}-{(x_i^k)}^{p-1}|
        = |(p-1){(\tilde{x}_i^k)}^{p-2}(\epsilon^{k-1}-x_i^k)|.
    \end{align*}
   For any $i\in\mathcal{I}2^k$, there exists  $\tilde{x}_i^k$ between $x_i^k$ and $x_i^{k-1}$ such that
    \begin{align*}
        |{(x_i^{k-1})}^{p-1}-{(x_i^k)}^{p-1}|
        = |(p-1){(\tilde{x}_i^k)}^{p-2}(x_i^{k-1}-x_i^k)|.
    \end{align*}
    It then follows that
    \begin{align*}
     \alpha(\boldsymbol{x}^k, \lambda^k) 
     \leq & \lambda^kp(1-p) \sum_{i\in\mathcal{I}1^k} {(\tilde{x}_i^k)}^{p-2}|\epsilon^{k-1}-x_i^k||x_i^k| \\
     & +\lambda^kp(1-p) \sum_{i\in\mathcal{I}2^k} {(\tilde{x}_i^k)}^{p-2}|x_i^{k-1}-x_i^k||x_i^k| \\
     \leq & \lambda^kp(1-p)\sum_{i\in\mathcal{I}1^k} {(\tilde{x}_i^k)}^{p-2}(|\epsilon^{k-1}-x_i^{k-1}|+|x_i^{k-1}-x_i^k|)|x_i^k|\\
     & +\lambda^kp(1-p) \sum_{i\in\mathcal{I}2^k} {(\tilde{x}_i^k)}^{p-2}|x_i^{k-1}-x_i^k||x_i^k| \\
     \leq & \lambda^kp(1-p)\sum_{i\in\mathcal{I}1^k} {(\tilde{x}_i^k)}^{p-2}\epsilon^{k-1}|x_i^k|\\
     & +\lambda^kp(1-p) \sum_{i\in\mathcal{I}^k} {(\tilde{x}_i^k)}^{p-2}|x_i^{k-1}-x_i^k||x_i^k|\\
     \leq & \lambda^kp(1-p)\epsilon^{k-1}\Vert{[\tilde{\boldsymbol{x}}^k_{\mathcal{I}1^k}]}^{p-2}\Vert_1\Vert \boldsymbol{y}\Vert_\infty\\
     & +\lambda^kp(1-p)\Vert{[\tilde{\boldsymbol{x}}^k_{\mathcal{I}^k}]}^{p-2}\Vert_2\Vert{\boldsymbol{x}}^{k-1}_{\mathcal{I}^k}-{\boldsymbol{x}}^k_{\mathcal{I}^k}\Vert_2\Vert \boldsymbol{y}\Vert_\infty,
    \end{align*}
    where the last inequality is by Cauchy-Schwarz inequality.

    Since $\boldsymbol{x}^k$ is the optimal solution of ($P_{k-1}$), by Lemma \ref{lemma5} (ii) and the definition of $\gamma^{k-1}$, we have
    \begin{align*}
    \beta(\boldsymbol{x}^k)=&|\sum_{i=1}^n (x_i^k)^p-\gamma|\\
    =&\gamma^{k-1} + \sum_{i=1}^n \phi^{k-1}(x_i^{k-1})-\sum_{i\in\mathcal{A}_{\epsilon}^{k-1}} 
    p{(\epsilon^{k-1})}^{p-1}x_i^k \\ & - \sum_{i\in\mathcal{I}_{\epsilon}^{k-1}} 
    p{(x_i^{k-1})}^{p-1}x_i^k - \sum_{i=1}^n (x_i^k)^p\\
    =&\sum_{i=1}^n \phi^{k-1}(x_i^{k-1})+
    \sum_{i\in\mathcal{A}_{\epsilon}^{k-1}} 
    p{(\epsilon^{k-1})}^{p-1}(x_i^k-x_i^{k-1}) \\ 
    & + \sum_{i\in\mathcal{I}_{\epsilon}^{k-1}}  
    p{(x_i^{k-1})}^{p-1}(x_i^k-x_i^{k-1})
    - \sum_{i=1}^n (x_i^k)^p\\
    \leq&\sum_{i=1}^n p{(\epsilon^{k-1})}^{p-1}|x_i^k- x_i^{k-1}| 
    +\sum_{i\in\mathcal{I}^k}|\phi^{k-1}(x_i^{k-1})-{(x^k_i)}^p|\\
    &+\sum_{i\in\mathcal{A}^k}\phi^{k-1}(x_i^{k-1}).
    \end{align*}
    %For any $i\in\mathcal{I}2^k$, by the Lagrange's mean value theorem, there exists $\tilde{x}_i^k$ between $x_i^k$ and $x_i^{k-1}$ such that
    %\begin{align*}
     %   |{(x^{k-1}_i)}^p-{(x^k_i)}^p|
     %   = |p{(\tilde{x}_i^k)}^{p-1}(x_i^{k-1}-x_i^k)|.
    %\end{align*}
    For any $i\in\mathcal{I}3^k$, according to (\ref{eqsum}), there exists $\tilde{x}_i^k$ between $x_i^k$ and $x_i^{k-1}+\epsilon^{k-1}$ such that
    \begin{align*}
        |\phi^{k-1}(x_i^{k-1})-{(x^k_i)}^p|
        & \leq |{{(x_i^{k-1}}+\epsilon^{k-1})}^{p} - {(x_i^k)}^{p}|\\
        & = |p{(\tilde{x}_i^k)}^{p-1}(x_i^{k-1}+\epsilon^{k-1}-x_i^k)|.
    \end{align*}
    For any $i\in\mathcal{I}4^k$, according to (\ref{eqsum}), there exists $\tilde{x}_i^k$ between $x_i^k$ and $x_i^{k-1}$ such that
    \begin{align*}
        |\phi^{k-1}(x_i^{k-1})-{(x^k_i)}^p|
        & \leq |{(x_i^{k-1}})^{p} - {(x_i^k)}^{p}|\\
        & = |p{(\tilde{x}_i^k)}^{p-1}(x_i^{k-1}-x_i^k)|.
    \end{align*}
    It then follows that
    \begin{align*}
    \beta(\boldsymbol{x}^k)\leq & p({\epsilon^{k-1}}){}^{p-1}\Vert \boldsymbol{x}^k- {\boldsymbol{x}}^{k-1}\Vert_1 + p\sum_{i\in\mathcal{I}3^k} {(\tilde{x}_i^k)}^{p-1}(|x_i^{k-1}-x_i^k|+ \\
    & |\epsilon^{k-1}|) + p\sum_{i\in\mathcal{I}4^k} {(\tilde{x}_i^k)}^{p-1}|x_i^{k-1}-x_i^k|
    + \sum_{i\in\mathcal{A}^k}(x^{k-1}_i+\epsilon^{k-1})^{p}\\ %eq1
    \leq & p({\epsilon^{k-1}}){}^{p-1}\Vert \boldsymbol{x}^k- {\boldsymbol{x}}^{k-1}\Vert_1 + p\sum_{i\in\mathcal{I}3^k} {(\tilde{x}_i^k)}^{p-1}\epsilon^{k-1} + \\
    & p\sum_{i\in\mathcal{I}^k} {(\tilde{x}_i^k)}^{p-1}|x_i^{k-1}-x_i^k|
    + \sum_{i\in\mathcal{A}^k}({x^{k-1}_i}+{\epsilon^{k-1}})({\epsilon^{k-1}})^{p-1}\\ %eq2
    \leq & p({\epsilon^{k-1}}){}^{p-1}\Vert \boldsymbol{x}^k- {\boldsymbol{x}}^{k-1}\Vert_1 + p\epsilon^{k-1}\Vert{[\tilde{\boldsymbol{x}}_{\mathcal{I}3^k}^k]}^{p-1}\Vert_1 \\ & + p\Vert {[\tilde{\boldsymbol{x}}_{\mathcal{I}^k}^k]}^{p-1}\Vert_2\Vert{\boldsymbol{x}}_{\mathcal{I}^k}^{k-1}-\boldsymbol{x}_{\mathcal{I}^k}^k\Vert_2
    +({\epsilon^{k-1}}){}^{p-1} \Vert {{\boldsymbol{x}}^{k-1}_{\mathcal{A}^k}}\Vert_1\\ & + ({\epsilon^{k-1}})^{p},
    \end{align*}
where the first and third inequalities are by Cauchy-Schwarz inequality.
\end{proof}

The following theorem presents the main result regarding the global convergence of the proposed algorithm.\\

\begin{theorem} %% theorem 2
Assume 0 is not the cluster point of $\{\lambda^k\}_\mathcal{U}$. Let $\{(\boldsymbol{x}^k, \lambda^k)\}$ be the sequence produced by Algorithm \ref{lloyd}. For any cluster point $(\boldsymbol{x}^*, \lambda^*)$ of $\{(\boldsymbol{x}^k, \lambda^k)\}_{\mathcal{U}}$,  $(\boldsymbol{x}^*, \lambda^*)$ is the first-order stationary point for problem (\ref{eq4}). 
\end{theorem}

\begin{proof}
Given $\{(\boldsymbol{x}^k, \lambda^k)\}_{\mathcal{U}}$, there exists $\mathcal{S}\subseteq\mathcal{U}$ such that $\lim_{{k\rightarrow +\infty}, k\in\mathcal{S}} (\boldsymbol{x}^k, \lambda^k) = (\boldsymbol{x}^*, \lambda^*)$.  { Therefore}, using the upper bounds established in Lemma \ref{lemma9}, we have 
$\alpha(\boldsymbol{x}^*, \lambda^*)= \lim_{{k\rightarrow +\infty}, k\in\mathcal{S}}  \alpha(\boldsymbol{x}^k, \lambda^k)=0$ and $\beta(\boldsymbol{x}^*)= \lim_{{k\rightarrow +\infty}, k\in\mathcal{S}}  \beta(\boldsymbol{x}^k)=0$ by Lemma \ref{lemma6}, Lemma \ref{lemma7}, and Lemma \ref{lemma8}. Hence, $\{(\boldsymbol{x}^*, \lambda^*)\}$ is stationary for problem (\ref{eq4}).  Overall, for any convergent subsequence $\mathcal{S}\subseteq\mathcal{U}$, we have $\lim_{{k\rightarrow +\infty}, k\in\mathcal{S}}  \alpha(\boldsymbol{x}^k, \lambda^k)=0$ and $\lim_{{k\rightarrow +\infty}, k\in\mathcal{S}}  \beta(\boldsymbol{x}^k)=0$. { Therefore}, both $\alpha(\boldsymbol{x}^k, \lambda^k)$ and $\beta(\boldsymbol{x}^k)$ converge to 0 on $\mathcal{U}$.
\end{proof}

In this work, we demonstrate that when 0 is not a cluster point of $\{\lambda^k\}_\mathcal{U}$, the sequence generated by the proposed algorithm always converges to a stationary point for (\ref{eq4}). As established in \citep{yang2022towards}, determining whether the original problem (\ref{eq4}) has a first-order stationary point $(\boldsymbol{x}^*, \lambda^*)$ with $\lambda^*=0$ is NP-hard. This complexity extends to the case where 0 is a cluster point of $\{\lambda^k\}_\mathcal{U}$, making it NP-hard to verify whether any point is  stationary point for (\ref{eq4}).

\section{Numerical experiments}

In this section, we conduct numerical experiments to compare the performance of the proposed algorithm with the existing IRBP algorithm that has demonstrated superiority over current state-of-the-art methods.We assess the performance and robustness of the proposed algorithm through tests on both synthetic and real-life data. The numerical experiments were implemented in Python using the NumPy package. All computations were performed on a MacBook Pro running Mac OS Monterey 12.0.1 with 8GB memory.

\subsection{Illustrative examples}
\begin{figure}[h]
    \centering
    \includegraphics[width=1\columnwidth]{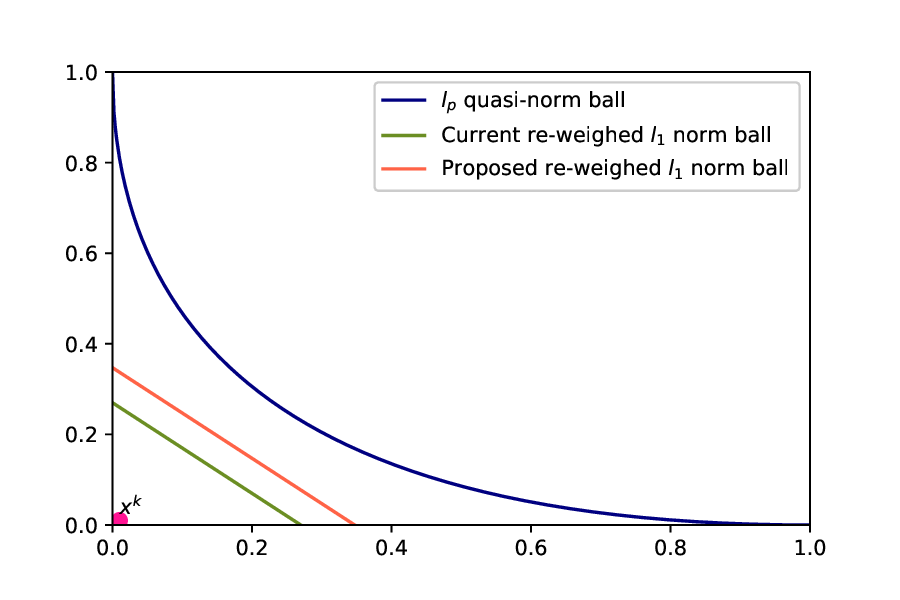}
    \caption{Illustrative example of the weighted $\ell_1$ norm balls produced by both IRBP algorithms ($n=2$, $p=0.5$, $\gamma=1$, $\boldsymbol{x}^k=(0.01,0.01)^\mathrm{T}$)}
    \label{fig1}
\end{figure}

To illustrate the efficacy of our approach, we present a two-dimensional illustrative example. Fig. \ref{fig1} depicts the feasible regions of the subproblem (\ref{Pk}) produced by two algorithms for the same problem instance with respect to the same iterate. It is observed that the two weighted $\ell_1$ norm balls are both contained within the original $\ell_p$ quasi-norm ball, but the one produced by the proposed algorithm encompasses the the other one. As a result, the proposed algorithm can achieve a more substantial improvement in the objective value in the next iteration. This enhancement can be attributed to the underlying approximate model $\sum_{i=1}^n \phi^k_i(x_i)$ in the proposed algorithm which provides a more refined localized approximation to the original $\ell_p$ quasi-norm ball. 

% as the iteration advances, $\epsilon^k$ is gradually reduced toward zero for a more precise approximation of the $\ell_p$ quasi-norm

%A noteworthy comparison between the existing IRBP algorithm and the proposed algorithm is that in the existing IRBP algorithm which employ the approximate model $\sum_{i=1}^n (|x_i|+\epsilon^k)^{p}$, for any feasible current iterate $\boldsymbol{x}^k$, the weight component corresponding to the zero component $x_i^k$ is computed as $p\epsilon^{p-1}$, and decreases accordingly as $x_i^k$ moves away from zero. In the proposed approximate model, for the iterate components whose value are near the origin within the interval [0, $u_{\epsilon^k}$], the corresponding weight component remains constant at $p\epsilon^{p-1}$, which is anticipated to more effectively promote sparsity by imposing a more substantial penalty for deviations from zero.

\begin{figure}[h]
    \centering
    \includegraphics[width=1\columnwidth]{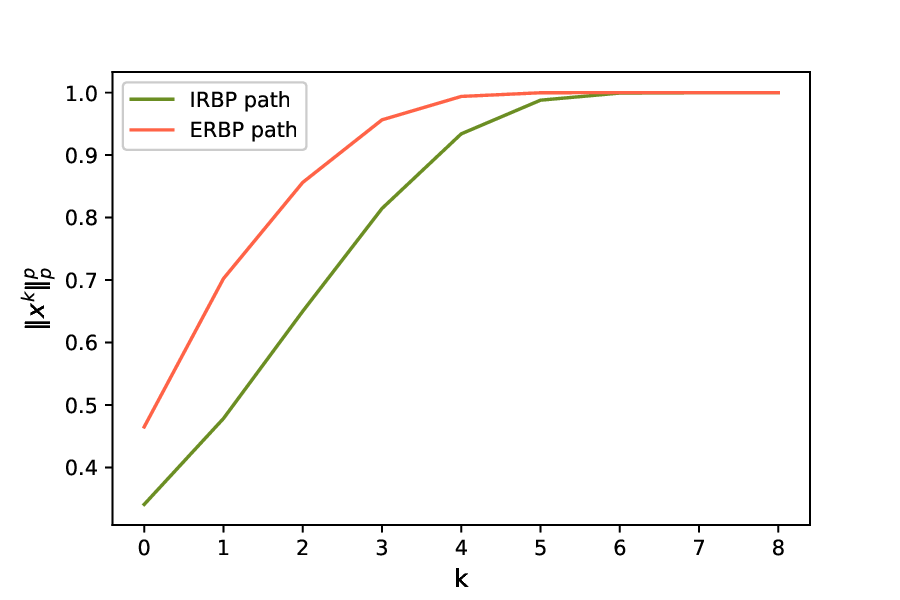}
    \caption{Illustrative example of the $\ell_p$ quasi-norm of the iterates produced by the IRBP and ERBP algorithms ($n=4$, $p=0.5$, $\gamma=1$, $\boldsymbol{y}=(0.18, 1.88, 0.20, 0.64)^\mathrm{T}$, $\boldsymbol{x}^0=(0,0,0,0)^\mathrm{T}$, $\epsilon=0.025$)}
    \label{fig2}
\end{figure}

For clarity, we will refer to the newly proposed algorithm as the Enhanced Iterative $\ell_1$ Norm Ball Projection Algorithm (ERBP). Figure \ref{fig2} compares the $\ell_p$ quasi-norm of the iterates $\boldsymbol{x}^k$ produced by the IRBP and ERBP algorithms as they progress. This comparison illustrates the rate at which each iterative sequence approaches the boundary of the original $\ell_p$ quasi-norm ball. The steeper slope observed in the ERBP path indicates that the ERBP algorithm is expected to require fewer iterations to approach the boundary of the quasi-norm ball than the IRBP algorithm.The improvement in the iteration path highlights the potential of the proposed method in achieving faster convergence and enhanced computational efficiency.

%This improvement in iteration path highlights the enhanced convergence characteristics of ERBP, making it a promising advancement in iterative optimization algorithms.
%In summary, the iteration path analysis underscores the superiority of ERBP, showcasing not only accelerated iteration time but also a more favorable iteration slope, collectively contributing to the algorithm's enhanced convergence performance.

\subsection{Experiments on synthetic data}
\label{section42}
To directly compare the effiicency and effecitveness of the proposed algorithm with the existing method, we first conduct projection experiments using synthetic data. We generated $\boldsymbol{y}$ with entries randomly sampled from a normal distribution with a mean of  $\mu$ and a standard deviation of one. Starting with $\mu=\gamma/n$, we incremented $\mu$ by $\gamma/n$ until obtaining a signal  $\boldsymbol{y}$ that lies outside the $\ell_p$ quasi-norm ball. We experimented with $p$ values of 0.4 and 0.6, reflecting different shapes of the quasi-norm ball. 

Regarding the input parameters of the algorithms, we initialized $\boldsymbol{x}^0$ as $\boldsymbol{0}$ and set $\epsilon^0$ to $0.4{(\gamma/n)}^{1/p}$ to ensure that $\phi^0(0)<\gamma/n$. During the $k$th iteration, we updated the value of $\delta$ using a strategy similar to those employed in the literature \citep{burke2015, yang2022towards}, as defined below:
\begin{equation}
\delta^{k+1} = \text{max}(1\text{e-}6, \text{min}(\beta(\boldsymbol{x}^k),1/\sqrt{k})^{1/p}).
\end{equation}
For both algorithms, we evaluated their performance using a stopping criterion
\begin{equation}
\alpha(\boldsymbol{x}^k,\epsilon^k)\leq \upsilon,\quad\beta(\boldsymbol{x}^k)\leq \upsilon,
\end{equation}
or if the number of iterations exceeded 1,000. We solved the weighted $\ell_1$ norm ball projection subproblems using the sort-based algorithm \citep{Konstantinos2010} with a computational complexity of order $O(n)$ for multiplications/additions and $O(n\log n)$ for sorting operations.

For each combination of input parameters $n$, $\gamma$, and $p$, we generated 20 signals to be projected, and ran the IRBP and ERBP algorithm independently with both $\upsilon$ values of \(1\text{e-}4\) and \(1\text{e-}8\) to thoroughly compare the convergence performance under different stopping criteria. 
In the first experiment, we examined the performance of two algorithms as the problem size increased. Specifically, we projected random vectors of varying sizes onto the $\ell_p$ quasi-norm ball with fixed radius. Table \ref{tab1} presents the problem instances along with the average CPU time (in seconds) and the number of iterations for both algorithms. { To faciliatete comparison, the method with better performance is highlight in bold.} Generally, the ERBP algorithm demonstrates lower computation time and requires fewer iterations than the IRBP algorithm. As the problem dimension scales up, the clear advantage in computational time and iterations demonstrated by the ERBP alrogithm becomes more pronounced. For a dimension of 1,000,000 with coarser and finer tolerances, ERBP outperforms by achieving approximately an 18\% and 17\% reduction in computational time for $p=0.4$ and a 29\% and 26\% reduction for $p=0.6$, respectively.

\begin{table*}[h]
\caption{Comparison of algorithms IRBP and ERBP  for varying problem sizes ($r=8$)}
\centering
{\begin{tabular}{c|ccccc|cccc}
\toprule    
& & \multicolumn{4}{c}{$\upsilon=1\text{e-}4$} &  \multicolumn{4}{c}{$\upsilon=1\text{e-}8$} \\
\midrule
& \multirow{2}{*}{$n$} & \multicolumn{2}{c}{CPU time} & \multicolumn{2}{c}{Iterations}   & \multicolumn{2}{c}{CPU time} & \multicolumn{2}{c}{Iterations}\\
& & IRBP & ERBP & IRBP & ERBP & IRBP & ERBP & IRBP & ERBP \\
\midrule
\multirow{6}{*}{$p=0.4$} 
 & 10 & \textbf{0.0015} & \textbf{0.0015} & 13.7 & \textbf{11.7} & \textbf{0.0020} & \textbf{0.0020} & 20.6 & \textbf{17.7} \\
 & 100 & 0.0034 & \textbf{0.0030} & 21.4 & \textbf{19.6} & 0.0036 & \textbf{0.0033} & 28.9 & \textbf{26.8} \\
 & 1,000 & 0.0118 & \textbf{0.0103} & 32.5 & \textbf{31.0} & 0.0144 & \textbf{0.0121} & 40.7 & \textbf{39.3} \\
 & 10,000 & 0.1031 & \textbf{0.0855} & \textbf{23.3} & \textbf{23.3} & 0.1122 & \textbf{0.0987} & \textbf{25.8} & \textbf{25.8} \\
 & 100,000 & 1.2677 & \textbf{1.1211} & 27.3 & \textbf{27.0} & 1.4187 & \textbf{1.2320} & \textbf{29.9} & \textbf{29.4} \\
 & 1,000,000 & 14.381 & \textbf{11.8401} & 33.3 & \textbf{32.5} & 16.5237 & \textbf{13.7503} & 38.9 & \textbf{37.2} \\
\midrule
\multirow{6}{*}{$p=0.6$} 
& 10 &  \textbf{0.0012} & \textbf{0.0012} & 11.9 & \textbf{10.2} & 0.0020 &  \textbf{0.0019} & 15.2 & \textbf{14.9} \\
 & 100 & 0.0014 & \textbf{0.0013} & 13.4 & \textbf{11.2} &  \textbf{0.0018} &  \textbf{0.0018} & 15.9 & \textbf{14.2} \\
 & 1,000 & 0.0049 & \textbf{0.0040} & 15.2 & \textbf{13.6} & 0.0057 &  \textbf{0.0047} & 16.6 & \textbf{15.4} \\
 & 10,000 & 0.0672 & \textbf{0.0476} & 15.6 & \textbf{13.4} & 0.0736 & \ \textbf{0.0555} & 16.4 & \textbf{14.6} \\
 & 100,000 & 0.7559 & \textbf{0.6231} & 16.2 & \textbf{15.2} & 0.8506 &  \textbf{0.7094} & 18.2 & \textbf{17.0} \\
 & 1,000,000 & 7.0539 & \textbf{4.9965} & 16.8 & \textbf{15.5} & 8.9969 &  \textbf{6.6827} & 20.5 & \textbf{18.3} \\
\bottomrule %   
\end{tabular}}
\label{tab1}
\end{table*}

In the second experiment, we focused on comparing the performances of two algorithms as the radius of the $\ell_p$ quasi-norm ball varies. Table \ref{tab2} illustrates that both algorithms are sensitive to the radius, with larger radii leading to increased computational time and iterations. Overall, the ERBP algorithm demonstrates lower computational time and generally requires fewer iterations compared to the IRBP algorithm across all tested radii. Specifically, with a radius of 128 and $p=0.4$, ERBP exhibits superior convergence performance by achieving approximately a 29\% and 26\% reduction in computational time for coarser and finer tolerances, respectively. For $p=0.6$ at the same radius, ERBP achieves a 16\% reduction for $p=0.6$ for both tolerances.

%This highlights that while ERBP shows an advantage in rapidly approaching from the interior to the boundary of the quasi-norm ball, it may experience slowdowns when iterating near the boundary and striving to meet stringent stopping criteria.

\begin{table*}[htb]
\caption{Comparison of algorithms IRBP and ERBP for varying $\ell_p$ quasi-norm ball radius ($n=10,000$)}
\centering
{\begin{tabular}{c|ccccc|cccc}
\toprule    
& & \multicolumn{4}{c}{$\upsilon=1\text{e-}4$} &  \multicolumn{4}{c}{$\upsilon=1\text{e-}8$} \\
\midrule
& \multirow{2}{*}{$\gamma$} & \multicolumn{2}{c}{CPU time} & \multicolumn{2}{c}{Iterations}   & \multicolumn{2}{c}{CPU time} & \multicolumn{2}{c}{Iterations}\\
& & IRBP & ERBP & IRBP & ERBP & IRBP & ERBP & IRBP & ERBP \\
\midrule
\multirow{8}{*}{$p=0.4$} 
 & 1 & 0.0038 & \textbf{0.0034} & \textbf{10.0} & \textbf{10.0} & 0.0043 & \textbf{0.0033} & 12.5 & \textbf{11.3} \\
 & 2 & 0.0062 & \textbf{0.0053} & 15.3 & \textbf{14.7} & 0.0086 & \textbf{0.0076} & 26.5 & \textbf{25.8} \\
 & 4 & 0.0064 & \textbf{0.0057} & 18.6 & \textbf{18.4} & 0.0095 & \textbf{0.0086} & 29.5 & \textbf{29.1} \\
 & 8 & 0.0103 & \textbf{0.0090} & 31.3 & \textbf{30.0} & 0.0121 & \textbf{0.0107} & 36.7 & \textbf{35.3} \\
 & 16 & 0.0128 & \textbf{0.0116} & 38.0 & \textbf{37.4} & 0.0133 & \textbf{0.0120} & 40.1 & \textbf{39.0} \\
 & 32 & 0.0216 & \textbf{0.0177} & 61.9 & \textbf{54.6} & 0.0231 & \textbf{0.0195} & 65.7 & \textbf{59.7} \\
 & 64 & 0.0338 & \textbf{0.0274} & 88.0 & \textbf{75.3} & 0.0342 & \textbf{0.0280} & 89.4 & \textbf{76.9} \\
 & 128 & 0.0505 & \textbf{0.0361} & 108.9 & \textbf{79.3} & 0.0515 & \textbf{0.0381} & 111.9 & \textbf{84.2} \\
\midrule
\multirow{8}{*}{$p=0.6$} 
 & 1 & 0.0032 & \textbf{0.0022} & 10.0 & \textbf{8.1} & 0.0035 & \textbf{0.0028} & 11.0 & \textbf{10.0} \\
 & 2 & 0.0043 & \textbf{0.0033} & 13.5 & \textbf{11.8} & 0.0095 & \textbf{0.0085} & 29.8 & \textbf{28.7} \\
 & 4 & 0.0056 & \textbf{0.0041} & 15.7 & \textbf{14.3} & 0.0056 & \textbf{0.0046} & 17.5 & \textbf{16.0} \\
 & 8 & 0.0049 & \textbf{0.0040} & 15.3 & \textbf{13.9} & 0.0065 & \textbf{0.0049} & 18.5 & \textbf{16.5} \\
 & 16 & 0.0051 & \textbf{0.0040} & 15.7 & \textbf{13.6} & 0.0057 & \textbf{0.0045} & 17.4 & \textbf{15.0} \\
 & 32 & 0.0056 & \textbf{0.0043} & 16.5 & \textbf{13.7} & 0.0061 & \textbf{0.0050} & 17.7 & \textbf{15.7} \\
 & 64 & 0.0064 & \textbf{0.0049} & 16.1 & \textbf{13.8} & 0.0066 & \textbf{0.0055} & 17.7 & \textbf{15.5} \\
 & 128 & 0.0074 & \textbf{0.0062} & 16.7 & \textbf{14.2} & 0.0081 & \textbf{0.0069} & 18.2 & \textbf{15.6} \\
\bottomrule %   
\end{tabular}}
\label{tab2}
\end{table*}

These results show that ERBP consistently outperforms IRBP across various problem sizes and radii, significantly enhancing the iteration path towards the boundary of the quasi-norm ball. This improvement is particularly pronounced in scenarios involving higher dimensions and larger radii. These findings are notable especially given prior research that has established the faster performance of the IRBP algorithm over other state-of-the-art methods \citep{yang2022towards}. With its superior computational efficiency, ERBP emerges as a promising method for computing the projection of a signal onto the $\ell_p$ quasi-norm ball, especially when faster convergence is needed. 

\subsection{Experiments on real-life data}
Given the practical relevance of computing projections onto the non-convex $\ell_p$ quasi-norm ball in domains such as signal processing and machine learning, we extend our empirical evaluation to address real-world problems in the context of image reconstruction. The approach is motivated by the sparsity of natural images in wavelet domains, which enables compressive imaging with a reduced number of measurements. Our experiments are conducted on eight 256$\times$256 grayscale images from the Set12 dataset \cite{7839189}. We employ the discrete wavelet transform (DWT) to construct sparse bases from the original images. Each of the 256$\times$256 matrix of wavelet coefficients is treated as 256 separate columns. For each column denoted as $\hat{\boldsymbol{x}}$, we compute $\boldsymbol{y} = A\hat{\boldsymbol{x}}$, where $A\in\mathbb{R}^{m\times 256}$ is a random measurement matrix with entries drawn from a standard normal distribution and the number of measurements $m$ is set to 200. We then solve the $\ell_p$-constrained least squares problem (\ref{eq3}) to obtain the estimated wavelet coefficients. After estimating across all columns, we perform the inverse transformation to recover the image. 

Since problem (\ref{eq3}) for $0<p<1$ is non-convex and challenging to solve globally, we approximate solutions using the projected gradient descent (PGD) method \cite{BAHMANI2013366}. The iterative update is based on projecting $\boldsymbol{x}^k + \eta^kA^T(\boldsymbol{y}-A\boldsymbol{x}^k)$ onto the $\ell_p$ quasi-norm ball of radius $\Vert\hat{\boldsymbol{x}}^k\Vert_p^p$, where $\eta^k > 0$ denotes the step size at iteration $k$. We integrated both IRBP and ERBP algorithms into the PGD framework for the projection step, and evaluated $p$ values of 0.4 and p=0.6. The experimental setup generally follows the outline in Section \ref{section42}. We report the average CPU time (in seconds) and Peak Signal-to-Noise Ratio (PSNR) over 20 image reconstruction trials in Table \ref{tab3}, where higher PSNR values indicate better quality of the the reconstructed image. Results suggest that the proposed ERBP algorithm in genreral reduces computational time across different datasets and $p$ values while maintaining comparable reconstruction quality, which highlights the practical value of this proposed approach in enhancing image recovery efficiency.

\subsection{Robust Weight Computation}
In numerical experiments, to address potential overflow issues, we introduced a safeguard by adding a small value, specifically $1\text{e-}12$, to the base of $(\cdot)^{p-1}$ to prevent the floating-point errors. Alternatively, when using an even smaller value like $1\text{e-}24$ for the base of $(\cdot)^{p-1}$, the ERBP algorithm successfully handled all problem instances and maintained consistent computational performance, but an increase in the computational time and iterations was noticeable for the IRBP algorithm, particularly as problem dimension or quasi-norm ball radius increased. In several instances, IRBP even failed to converge within the allotted 1,000 iterations due to adjustments in weight calculation, potentially leading to trapping during iteration despite its theoretical global convergence. This underscores the practical computational robustness for the proposed ERBP algorithm, particularly advantageous in very high-dimensional problems.

\begin{table*}[h]
\caption{Comparison of algorithms IRBP and ERBP for the image reconstruction tests}
\centering
{\begin{tabular}{c|cccc|cccc}
\toprule    
& \multicolumn{4}{c}{$p=0.4$} &  \multicolumn{4}{c}{$p=0.6$} \\
\midrule
\multirow{2}{*}{Dataset} & \multicolumn{2}{c}{CPU time} & \multicolumn{2}{c}{PNSR}   & \multicolumn{2}{c}{CPU time} & \multicolumn{2}{c}{PNSR}\\
& IRBP & ERBP & IRBP & ERBP & IRBP & ERBP & IRBP & ERBP \\
\midrule
lena & 262.02 & \textbf{177.19} & \textbf{30.53} & 30.37 & 846.59 & \textbf{721.88} & \textbf{30.55} & 30.33 \\ 
butterfly & 303.22 & \textbf{199.92} & 26.90 & \textbf{27.59} & 1009.76 & \textbf{815.79} & 27.13 & \textbf{27.17} \\ 
cameraman & 231.17 & \textbf{153.40} & 33.58 & \textbf{34.13} & 684.51  & \textbf{568.44} & 34.18 & \textbf{34.35} \\ 
peppers & 260.78 & \textbf{157.54} & 31.09 & \textbf{31.24} & 714.37  & \textbf{639.55} & 30.41 & \textbf{31.78} \\ 
house & 138.03 & \textbf{90.90} & \textbf{35.67} & 35.66 & 356.82  & \textbf{316.58} & \textbf{35.97} & 35.85 \\ 
jetplane & 284.25 & \textbf{186.75} & 29.69 & \textbf{29.83} & \textbf{966.46}  & 1149.21 & 27.84 & \textbf{29.03} \\ 
sea star & 340.03 & \textbf{266.94} & 25.30 & \textbf{26.08} & 1862.01  & \textbf{1372.39} & 25.25 & \textbf{25.61} \\ 
parrot & 275.17 & \textbf{181.53} & \textbf{29.38} & 29.34 & 2178.49  & \textbf{768.90} & \textbf{29.92} & 29.89 \\ 
\bottomrule 
\end{tabular}}
\label{tab3}
\end{table*}

\section{Summary}

The problem of projecting onto the $\ell_p$ quasi-norm ball is fundamental in various machine learning applications, such as signal denoising and data compression. However,  effective algorithms for its solution remain limited. In this paper, we introduce a novel variant of a state-of-the-art iteratively reweighted algorithm, specifically designed to solve the $\ell_p$ quasi-norm ball projection through iterative projections onto a sequence of weighted $\ell_1$ norm balls. A key challenge is the non-Lipschitz continuity of the $\ell_p$ quasi-norm function. While previous methods have addressed this by smoothing the function over the entire domain, we propose a localized approximation that focuses specifically on the neighborhood around the origin, ensuring Lipschitz continuity is preserved throughout the entire domain. This strategy enables us to achieve an $\epsilon$-approximation of the $\ell_p$ quasi-norm, which not only improves the quality of approximation but also preserves the concave structure inherent in the original $\ell_p$ quasi-norm function. By leveraging these properties, we enhance the existing iterative reweighted $\ell_1$ norm ball projection method by constructing a sequence of tighter subproblems that effectively approach the boundary of the non-convex quasi-norm ball from the interior. The global convergence properties of the proposed algorithm have been rigorously established. In direct comparison with existing approaches on projection experiment using random data, the proposed algorithm consistently demonstrates improved computational efficiency in terms of both computational time and number of iterations. Moreover, in image reconstruction experiments with real-life data, the proposed algorithm consistently enhances image recovery efficiency while maintaining high reconstruction quality. These results strongly support the proposed algorithm as a promising, efficient method for computing projections onto the $\ell_p$ quasi-norm ball with substantial practical value. \\ 

\noindent\textbf{Competing interests}
The authors have no known competing financial interests or personal relationships that could have appeared to influence the work reported in this paper.\\

\noindent\textbf{Authors contribution statement} Conceptualization, Methodology: QA; Original draft preparation: QA, JW. Data analysis, Writing - review and editing: ZN, NZ. All authors approved the final manuscript.\\

%\noindent\textbf{Ethical and informed consent for data used} 
%This article does not contain any studies with human participants or animals performed by any of the authors.

\noindent\textbf{Funding} 
The research is funded by the Innovation Fund of the Engineering Research Center of Integration and Application of Digital Learning Technology, Ministry of Education (1421014).\\

\noindent\textbf{Data availability and access} 
All relevant data are available at https://github.com/nivalann721/lpquasinormprojection.\\

%\bibliographystyle{cas-model2-names}
%\bibliography{main}

\end{document}